\documentclass[runningheads]{llncs}
\usepackage[T1]{fontenc}
\usepackage{graphicx}
\usepackage{amssymb}
\usepackage{amsmath}
\usepackage{subfigure}

\usepackage{amsbsy}
\usepackage{amsfonts}
\usepackage{mathbbol}

\DeclareMathOperator*{\argmin}{argmin}

\begin{document}
\title{
Some smooth divergences for $\ell_{1}-$approximations 
}
\titlerunning{Smooth divergences}
%
\author{
{Pierre Bertrand
\and
Wolfgang Stummer
}
}
\authorrunning{P. Bertrand and W. Stummer
}
%
\institute{
AMSE, Aix-Marseille Universit\'{e}, 5-9 Boulevard Maurice Bourdet, CS 50498, 13205 Marseille, France,
\and
Corresponding Author.
Department of Mathematics, Friedrich-Alexander-Universit{\"at} Erlangen--N\"{u}rnberg,
Cauerstrasse $11$, 91058 Erlangen, Germany,
\email{stummer@math.fau.de}. 
\\[-0.7cm]
\ 
}
\maketitle              
\begin{abstract}
For some smooth special case of generalized $\varphi-$divergences
as well as of new divergences (called scaled shift divergences),
we derive approximations of the omnipresent
(weighted) $\ell_{1}-$distance and (weighted) $\ell_{1}-$norm. 

\\[-1.0cm]
\ 
\end{abstract}


\enlargethispage{2.2cm}

\section{The generalized $\varphi-$divergence case}
\label{BerBroStuGSI25:sec1}

\vspace{-0.2cm} 
\noindent
It is well-known that a divergence is a real-valued function $D$ on
(a subset of) $\mathbb{R}^{K} \times \mathbb{R}^{K}$
which has the following two properties: (i) $D( \mathbf{Q}, \mathbf{P} ) \geq 0$ for 
$K-$dimensional vectors $\mathbf{Q}, \mathbf{P}$, 
and (ii) $D( \mathbf{Q}, \mathbf{P} ) = 0$ if and only if $\mathbf{Q} = \mathbf{P}$.
Since in general, $D( \mathbf{Q}, \mathbf{P} ) \ne D( \mathbf{P}  , \mathbf{Q} )$
and the triangle inequality is not satisfied,
$D( \mathbf{Q}, \mathbf{P} )$ can be interpreted as a \textit{directed distance}; 
accordingly, the divergence $D$ can be connected to geometric issues in various different ways,
see e.g. the detailed discussion in Section 1.5 of \cite{BerBroStuGSI25:Bro22b},
and \cite{BerBroStuGSI25:Roe17}.
Typically, a divergence $D$ is generated by some function $\varphi$. For the latter,
we require for the rest of the paper: 

\vspace{-0.1cm}

\begin{itemize}

\item
$\varphi: \, ]-\infty,
\infty[ \rightarrow [0,\infty]$ is lower semicontinuous and convex,
\ with $\varphi(1)=0$;

\item
the effective domain $dom(\varphi) := \{ t \in \mathbb{R} : \varphi(t) < \infty \}$
has interior $int(dom(\varphi))$
of the form $int(dom(\varphi)) = \, ]a,b[$ for
some $-\infty \leq a < 1 < b \leq \infty$;

\item
$t \mapsto \varphi(t)$ is strictly convex at $t=1$ 
(i.e. it is not identically zero   
in the open interval $]1-\varepsilon,1+\varepsilon[$
for any $\varepsilon > 0$,
cf. e.g. Liese \& Miescke \cite{BerBroStuGSI25:Lie08}).

\end{itemize}

\vspace{-0.1cm}
\noindent
Also, we set $\varphi(a) := \lim_{t \downarrow a} \varphi(t)$ 
and $\varphi(b) := \lim_{t \uparrow b} \varphi(t)$. 
For $\mathbf{P}:=\left( p_{1},\ldots,p_{K}\right) \in 
\mathbb{R}_{>0}^{K}
:= \{\mathbf{R}:= (r_{1},\ldots,r_{K}) \in \mathbb{R}^K: \,
r_{i} > 0 \ \text{for all} \ i=1,\ldots,K \}
$
and $\mathbf{Q} := (q_{1},\ldots,q_{K}) \in \mathbf{\Omega} \subset \mathbb{R}^{K}$, we define as directed distance the \textit{generalized $\varphi-$divergence} (generalized Csisz\'ar-Ali-Silvey-Morimoto divergence) 
\vspace{-0.2cm}
\begin{equation}
D_{\varphi}( \mathbf{Q}, \mathbf{P} ) := 
\sum_{k=1}^{K} p_{k} \cdot
\varphi \left( \frac{q_{k}}{p_{k}}\right); 
\label{BerBroStuGSI25:fo.genCASM}
\end{equation}

\vspace{-0.2cm}
\noindent
for a comprehensive technical
treatment, see e.g. 
\cite{BerBroStuGSI25:Bro19a}.
Comprehensive overviews on these important (generalized) $\varphi-$divergences
are given in 
e.g.~\cite{BerBroStuGSI25:Lie87},~\cite{BerBroStuGSI25:Vaj89},~\cite{BerBroStuGSI25:Bro22b},~\cite{BerBroStuGSI25:Bro23a},
and the references therein. Notice that the 
$\ell_{1}-$distance --- also called total variation distance ---
$D_{\varphi_{TV}}( \mathbf{Q}, \mathbf{P} ) := 
\sum_{k=1}^{K} p_{k} \cdot
\varphi_{TV} \negthinspace \left( \frac{q_{k}}{p_{k}}\right) = 
\sum_{k=1}^{K} | \, p_{k} - q_{k}\, |
$ with $\varphi_{TV}(t) := |t-1|$ is covered here.
Another interesting example is given as follows:
for any parameter-triple $\alpha,\beta,\widetilde{c} \in \, ]0,\infty[$
we choose $]a,b[ \,  :=  \,  ]-\infty, \infty \, [$  and

\begin{eqnarray}
\varphi_{\alpha,\beta,\widetilde{c}}(t) 
&:=& 
\begin{cases}
\widetilde{c} \cdot \alpha \cdot \Big\{
\sqrt{1 + \beta^{2} 
\cdot \Big(\frac{1-t}{\alpha}\Big)^2} - 1 
+ \log\frac{2 \cdot \Big(
\sqrt{1 + \beta^{2} 
\cdot \Big(\frac{1-t}{\alpha} \Big)^2} - 1
\Big)
}{
\beta^{2} \cdot \Big(\frac{1-t}{\alpha}\Big)^{2}
}  \Big\} \in \, ]0,\infty[,
\\
\hspace{6.5cm} \textrm{if } \  
t \in \, ]-\infty, 1 [ \, \cup \, ]1,\infty[ \, , 
\\
0, \hspace{6.2cm}  
\textrm{if } \  
t=1, 
\end{cases}
\label{BerBroStuGSI25:fo.genLap3ba.BS3}
\end{eqnarray}

\vspace{-0.3cm}
\noindent
(cf. Broniatowski \& Stummer \cite{BerBroStuGSI25:Bro23a}).
Notice that  
$\varphi_{\alpha,\beta,\widetilde{c}}(1) = 0$, 
$\varphi_{\alpha,\beta,\widetilde{c}}^{\prime}(1) = 0$,
$\varphi_{\alpha,\beta,\widetilde{c}}(-\infty) = \infty$
and $\varphi_{\alpha,\beta,\widetilde{c}}(\infty) = \infty$. 
Moreover,  
$\varphi_{\alpha,\beta,\widetilde{c}}^{\prime}(-\infty) = 
\varphi_{\alpha,\beta,\widetilde{c}}^{\prime}(a) =
-\widetilde{c} \cdot \beta$ and 
$\varphi_{\alpha,\beta,\widetilde{c}}^{\prime}(\infty) = 
\varphi_{\alpha,\beta,\widetilde{c}}^{\prime}(b) =
\widetilde{c} \cdot \beta$.  Furthermore, 
$\varphi_{\alpha,\beta,\widetilde{c}}(\cdot)$ is strictly convex and
smooth (i.e. of $C^{\infty}-$type), and
$\varphi_{\alpha,\beta,\widetilde{c}}(t) \leq \widetilde{c} \cdot \beta \cdot |t-1|$
with equality iff $t=1$.
From \eqref{BerBroStuGSI25:fo.genLap3ba.BS3},
we construct the generalized $\varphi-$divergence 
\vspace{-0.2cm}
\begin{eqnarray}
& & 
\hspace{-0.5cm}
D_{\varphi_{\alpha,\beta,\widetilde{c}}}(\mathbf{Q},\mathbf{P})
= \sum\limits_{k=1}^{K} p_{k} \cdot 
\varphi_{\alpha,\beta,\widetilde{c}}\Big(\frac{q_{k}}{p_{k}}\Big)
\nonumber \\[-0.1cm]
& & 
\hspace{-0.5cm}
= 
\begin{cases}
\widetilde{c} \cdot \alpha \cdot  \sum\limits_{k=1}^{K} 
p_{k} \cdot
\Big\{
\sqrt{1 + \beta^{2} 
\cdot \Big(\frac{1-\frac{q_{k}}{p_{k}}}{\alpha}\Big)^2} \, - \, 1 
+ \log\frac{2 \cdot \Big(
\sqrt{1 + \beta^{2} 
\cdot \Big(\frac{1-\frac{q_{k}}{p_{k}}}{\alpha} \Big)^2} \, - \, 1
\Big)
}{
\beta^{2} \cdot \Big(\frac{1-\frac{q_{k}}{p_{k}}}{\alpha}\Big)^{2}
}  \Big\} , 
\\ 
\hspace{6.6cm} \textrm{if } \  \mathbf{P} \in \mathbb{R}_{> 0}^{K},
\mathbf{Q} \in \mathbb{R}^{K}\backslash\{\mathbf{P}\} \, , 
\\
0, \hspace{6.25cm}  
\textrm{if } \  
\mathbf{Q} = \mathbf{P}. 
\end{cases}
\label{BerBroStuGSI25:fo.genLap4fa.BS3}
\end{eqnarray}

\vspace{-0.2cm}
\noindent
As a background, in \cite{BerBroStuGSI25:Bro23a}
we have shown that 
for any fixed $\mathbf{P} \in \mathbb{R}_{> 0}^{K}$ with 
$M_{\mathbf{P}} =\sum_{i=1}^{K}p_{i} \in \, ]0,\infty[$
there holds for $\varphi := \varphi_{\alpha,\beta,\widetilde{c}}$
the important condition
\vspace{-0.2cm} 
\begin{equation}
M_{\mathbf{P}} \cdot \varphi(t)
= 
\sup_{z \in \mathbb{R}} \Big( z\cdot t - \log \int_{\mathbb{R}} e^{zy} d\widetilde{\mathbb{\zeta}}(y) \Big),
\qquad t \in \mathbb{R},  
\nonumber
\end{equation}

\enlargethispage{0.5cm}

\vspace{-0.3cm} 
\noindent
for some probability distribution $\widetilde{\mathbb{\zeta}}$ on the real line
such that the function $z \mapsto MGF_{\widetilde{\mathbb{\zeta}}}(z) := \int_{\mathbb{R}} e^{zy} d\widetilde{\mathbb{\zeta}}(y)$ is finite
on some open interval containing zero.
Indeed, the corresponding distribution $\widetilde{\zeta}[\,\cdot \,] :=
\widetilde{\zeta}_{\alpha,\beta,\widetilde{c}}[\,\cdot \,]
:= \Pi[\widetilde{W}\in \cdot \,]$ 
is the comfortably simulable \textit{generalized 
Laplace distribution} 
of a random variable $\widetilde{W}:= 1 + \widetilde{Z}_{1} - \widetilde{Z}_{2}$, where
$\widetilde{Z}_{1}$ and $\widetilde{Z}_{2}$ are auxiliary random variables which are 
independent and identically
$GAM(M_{\mathbf{P}} \cdot \widetilde{c} \cdot \beta,M_{\mathbf{P}} \cdot \widetilde{c} \cdot \alpha)-$distributed.  
Accordingly, $\widetilde{W}$ has expectation $1$ and variance 
$2 \cdot \frac{(M_{\mathbf{P}} \cdot \widetilde{c} \cdot \alpha)^{2}}{M_{\mathbf{P}} \cdot \widetilde{c} \cdot \beta}
= 2 M_{\mathbf{P}} \cdot \widetilde{c} \cdot \frac{\alpha^{2}}{\beta}$.

\vspace{0.2cm}
\noindent
In the following, we show how the 
$\varphi-$divergence
\eqref{BerBroStuGSI25:fo.genLap4fa.BS3} can be employed 
to achieve smooth approximations of the $\ell_{1}-$distance as well as 
the $\ell_{1}-$norm:

\vspace{-0.1cm}
\begin{proposition}
\label{BerBroStuGSI25:prop.limit.ell1.BS3}
(a) For all $t \in \, ]-\infty,\infty[$, $\beta \in \, ]0,\infty[$ and $\widetilde{c} \in \, ]0,\infty[$
there holds
\begin{equation}
\lim_{\alpha \rightarrow 0_{+}} \varphi_{\alpha,\beta,\widetilde{c}}(t) \ = \ 
\widetilde{c} \cdot \beta \cdot |t-1| .
\nonumber
\end{equation}

\vspace{-0.2cm}
\noindent
(b) For all $t \in \, ]-\infty,\infty[$, $\alpha \in \, ]0,\infty[$ and $\beta \in \, ]0,\infty[$ 
there holds
\begin{equation}
\lim_{\frac{\alpha}{\beta} \rightarrow 0_{+}} \varphi_{\alpha,\beta,1/\beta}(t) \ = \ 
|t-1| .
\nonumber
\end{equation}

\vspace{-0.2cm}
\noindent
(c) For all $\beta \in \, ]0,\infty[$, $\widetilde{c} \in \, ]0,\infty[$,
$\mathbf{Q} \in \mathbb{R}^{K}$ and $\mathbf{P} \in \mathbb{R}_{>0}^{K}$
there holds
\vspace{-0.2cm}
\begin{equation}
\lim_{\alpha \rightarrow 0_{+}}
D_{\varphi_{\alpha,\beta,\widetilde{c}}}(\mathbf{Q},\mathbf{P}) \ = \ 
\widetilde{c} \cdot \beta \cdot \sum_{k=1}^{K} |q_{k} - p_{k}| 
\ = \ \widetilde{c} \cdot \beta \cdot ||\mathbf{Q}-\mathbf{P}||_{1} .
\nonumber
\end{equation}
\noindent
(d) For all $\alpha \in \, ]0,\infty[$, $\beta \in \, ]0,\infty[$, 
$\mathbf{Q} \in \mathbb{R}^{K}$ and $\mathbf{P} \in \mathbb{R}_{>0}^{K}$
there holds
\vspace{-0.3cm}
\begin{equation}
\lim_{\frac{\alpha}{\beta} \rightarrow 0_{+}}
D_{\varphi_{\alpha,\beta,1/\beta}}(\mathbf{Q},\mathbf{P}) \ = \ 
\sum_{k=1}^{K} |q_{k} - p_{k}| 
\ = \ ||\mathbf{Q}-\mathbf{P}||_{1} .
\label{BerBroStuGSI25:fo.limit.ell1.2.BS3sec}
\end{equation}

\newpage

\enlargethispage{0.5cm}

\vspace{-0.3cm}
\noindent
(e) For all $\alpha \in \, ]0,\infty[$, $\beta \in \, ]0,\infty[$, $\widetilde{c} \in\, ]0,\infty[$,
$\mathbf{Q} \in \mathbb{R}^{K}$ and all sequences $(\mathbf{P}_{m})_{m \in \mathbb{N}}$ in $\mathbb{R}_{>0}^{K}$
which tend (component-wise) to $\mathbf{0}$  (i.e. $\mathbf{P}_{m} \stackrel{m\rightarrow \infty}{\longrightarrow}
\mathbf{0}$) there holds
\vspace{-0.3cm}
\begin{equation}
\lim_{m \rightarrow \infty}
D_{\varphi_{\alpha,\beta,\widetilde{c}}}(\mathbf{Q},\mathbf{P}_{m}) \ = \ 
\widetilde{c} \cdot \beta \cdot \sum_{k=1}^{K} |q_{k}| 
\ = \ \widetilde{c} \cdot \beta \cdot ||\mathbf{Q}||_{1} ;
\label{BerBroStuGSI25:fo.limit.ell1.2.BS3a}
\end{equation}

\vspace{-0.2cm}
\noindent 
especially for 
$\mathbf{P}_{m} := \frac{1}{m} \cdot \mathbf{1}$ 
(i.e. each of the $K$ components has value $\frac{1}{m}$)
one has
\vspace{-0.2cm}
\begin{equation}
\lim_{m \rightarrow \infty}
D_{\varphi_{\alpha,\beta,\widetilde{c}}}(\mathbf{Q},\frac{1}{m} \cdot \mathbf{1}) \ = \ 
\widetilde{c} \cdot \beta \cdot \sum_{k=1}^{K} |q_{k}| 
\ = \ \widetilde{c} \cdot \beta \cdot ||\mathbf{Q}||_{1} .
\nonumber
\end{equation}

\end{proposition}

\vspace{-0.2cm}
\noindent
The proof of Proposition \ref{BerBroStuGSI25:prop.limit.ell1.BS3}
will be given in Section \ref{BerBroStuGSI25:sec:4} below.


\vspace{-0.4cm}

\section{The case of the 
scaled shift-divergence case}

\vspace{-0.1cm}
\noindent
Instead of the \textit{generalized $\varphi-$divergence} \eqref{BerBroStuGSI25:fo.genCASM},
let us now construct 
--- for $\mathbf{P} \in \mathbb{R}_{> 0}^{K}$, 
$\mathbf{Q} \in \mathbb{R}^{K}$,
$\mathbf{Q}^{\ast} \in \mathbb{R}^{K}$ 
and $\mathbf{\sigma} \in \mathbb{R}_{> 0}^{K}$
---
the new \textit{scaled shift divergence}
\vspace{-0.2cm}
\begin{equation}
D_{\varphi,\mathbf{P},\mathbf{\sigma}}^{new}(\mathbf{Q},\mathbf{Q}^{\ast}) \ := \ 
\sum\limits_{k=1}^{K} p_{k} \cdot
\varphi\negthinspace 
\left( \frac{q_{k}-q_{k}^{\ast}}{p_{k} \cdot \sigma_{k}} + 1 \right)
\nonumber
\end{equation}

\vspace{-0.3cm}
\noindent
where the divergence generator $\varphi$ has the general properties 
declared in the beginning of Section \ref{BerBroStuGSI25:sec1}.
Notice that $\mathbf{P}$ plays a different role --- namely that of a weight/scaling ---
than in \eqref{BerBroStuGSI25:fo.genCASM}. Clearly, there hold the divergence properties
$D_{\varphi,\mathbf{P},\mathbf{\sigma}}^{new}(\mathbf{Q},\mathbf{Q}^{\ast}) \geq 0$,
with equality if and only if $\mathbf{Q} = \mathbf{Q}^{\ast}$.
For the special choice $\varphi := \varphi_{\alpha,\beta,\widetilde{c}}$ 
(cf. \eqref{BerBroStuGSI25:fo.genLap3ba.BS3}) we end up with
\vspace{-0.3cm}
\begin{eqnarray}
& & 
\hspace{-0.2cm}
D_{\varphi_{\alpha,\beta,\widetilde{c}},\mathbf{P},\mathbf{\sigma}}^{new}(\mathbf{Q},\mathbf{Q}^{\ast}) \ := \ 
\sum\limits_{k=1}^{K} p_{k} \cdot
\varphi_{\alpha,\beta,\widetilde{c}} \negthinspace 
\left( \frac{q_{k}-q_{k}^{\ast}}{p_{k} \cdot \sigma_{k}} + 1 \right)
\nonumber \\[-0.1cm]
& & 
\hspace{-0.2cm}
= 
\begin{cases}
\widetilde{c} \cdot \alpha \cdot  \sum\limits_{k=1}^{K} 
p_{k} \cdot
\Big\{
\sqrt{1 + \frac{\beta^{2}}{\alpha^{2}} 
\cdot \Big(\frac{q_{k}-q_{k}^{\ast}}{p_{k} \cdot \sigma_{k}}\Big)^2} - 1 
+ \log\frac{2 \cdot \Big(
\sqrt{1 + \frac{\beta^{2}}{\alpha^{2}} 
\cdot \Big(\frac{q_{k}-q_{k}^{\ast}}{p_{k} \cdot \sigma_{k}}\Big)^2} - 1
\Big)
}{
\frac{\beta^{2}}{\alpha^{2}} 
\cdot \Big(\frac{q_{k}-q_{k}^{\ast}}{p_{k} \cdot \sigma_{k}}\Big)^2
}  \Big\}   
 \in \, ]0,\infty[,
\\
\hspace{8.6cm} \textrm{if } \  
\mathbf{Q} \in \mathbb{R}^{K}\backslash\{\mathbf{Q}^{\ast}\} \, , 
\\
0, \hspace{8.3cm}  
\textrm{if } \  
\mathbf{Q} = \mathbf{Q}^{\ast}. 
\end{cases}
\nonumber
\end{eqnarray}

\vspace{0.2cm}
\noindent
For this, we can deduce the following \textit{weighted} $\ell_{1}-$distance approximations:

\begin{proposition}
\label{BerBroStuGSI25:prop.newdiv.ell1.BS3}
(a) For all $\beta \in \, ]0,\infty[$, $\widetilde{c} \in \, ]0,\infty[$,
$\mathbf{Q} \in \mathbb{R}^{K}$, $\mathbf{Q}^{\ast} \in \mathbb{R}^{K}$, $\mathbf{P} \in \mathbb{R}_{>0}^{K}$
and $\mathbf{\sigma} \in \mathbb{R}_{>0}^{K}$
there holds
\vspace{-0.3cm}
\begin{equation}
\lim_{\alpha \rightarrow 0_{+}}
D_{\varphi_{\alpha,\beta,\widetilde{c}},\mathbf{P},\mathbf{\sigma}}^{new}(\mathbf{Q},\mathbf{Q}^{\ast})
\ = \ 
\widetilde{c} \cdot \beta \cdot \sum_{k=1}^{K} \frac{|q_{k} - q_{k}^{\ast}|}{\sigma_{k}}
\nonumber
\end{equation}

\vspace{-0.2cm}
\noindent
where the latter is (a multiple of) a weighted $\ell_{1}-$distance between $\mathbf{Q}$ and $\mathbf{Q}^{\ast}$.\\
(b) For all $\alpha \in \, ]0,\infty[$, $\beta \in \, ]0,\infty[$, 
$\mathbf{Q} \in \mathbb{R}^{K}$, $\mathbf{Q}^{\ast} \in \mathbb{R}^{K}$, $\mathbf{P} \in \mathbb{R}_{>0}^{K}$,
$\mathbf{\sigma} \in \mathbb{R}_{>0}^{K}$ we get
\vspace{-0.2cm}
\begin{equation}
\lim_{\frac{\alpha}{\beta} \rightarrow 0_{+}}
D_{\varphi_{\alpha,\beta,1/\beta},\mathbf{P},\mathbf{\sigma}}^{new}(\mathbf{Q},\mathbf{Q}^{\ast})
\ = \ 
\sum_{k=1}^{K} \frac{|q_{k} - q_{k}^{\ast}|}{\sigma_{k}} .
\label{BerBroStuGSI25:prop.newdiv.ell1.2.BS3sec}
\end{equation}

\newpage

\vspace{-0.2cm}
\noindent
(c) For all $\alpha \in \, ]0,\infty[$, $\beta \in \, ]0,\infty[$, $\widetilde{c} \in\, ]0,\infty[$,
$\mathbf{Q} \in \mathbb{R}^{K}$, $\mathbf{Q}^{\ast} \in \mathbb{R}^{K}$, $\mathbf{\sigma} \in \mathbb{R}_{>0}^{K}$ 
and all sequences $(\mathbf{P}_{m})_{m \in \mathbb{N}}$ in $\mathbb{R}_{>0}^{K}$
which tend to $\mathbf{0}$  (i.e. $\mathbf{P}_{m} \stackrel{m\rightarrow \infty}{\longrightarrow}
\mathbf{0}$) there holds
\vspace{-0.2cm}
\begin{equation}
\lim_{m \rightarrow \infty}
D_{\varphi_{\alpha,\beta,\widetilde{c}},\mathbf{P}_{m},\mathbf{\sigma}}^{new}(\mathbf{Q},\mathbf{Q}^{\ast})
\ = \ 
\widetilde{c} \cdot \beta \cdot \sum_{k=1}^{K} \frac{|q_{k} - q_{k}^{\ast}|}{\sigma_{k}} ;
\label{BerBroStuGSI25:prop.newdiv.ell1.2.BS3a}
\end{equation}

\vspace{-0.2cm}
\noindent
in particular, for $\mathbf{P}_{m} := \frac{1}{m} \cdot \mathbf{1}$ 
there holds
\vspace{-0.2cm}
\begin{equation}
\lim_{m \rightarrow \infty}
D_{\varphi_{\alpha,\beta,\widetilde{c}},\frac{1}{m} \cdot \mathbf{1},\mathbf{\sigma}}^{new}(\mathbf{Q},\mathbf{Q}^{\ast})
\ = \ 
\widetilde{c} \cdot \beta \cdot \sum_{k=1}^{K} \frac{|q_{k} - q_{k}^{\ast}|}{\sigma_{k}} .
\nonumber
\end{equation}

\end{proposition}

\vspace{-0.2cm}
\noindent
For the special choice $\mathbf{Q}^{\ast} = \mathbf{0}$ in (a),(b),(c)
we immediately obtain the corresponding limit assertions for the
\textit{weighted $\ell_{1}-$norm}. 
The proof of Proposition \ref{BerBroStuGSI25:prop.newdiv.ell1.BS3}
will be given in Section \ref{BerBroStuGSI25:sec:4} below.

\vspace{-0.2cm}


\section{Some visualizations}

\vspace{-0.1cm}
\noindent 
In the following, we visualize some of the above convergences
in a concrete context, say, at a LASSO minimizer point
(cf. Tibshirani \cite{BerBroStuGSI25:Tib96}) 
\vspace{-0.2cm}
\begin{equation}
\widehat{\mathbf{Q}} \in \argmin_{\mathbf{Q} \in \mathbb{R}^{K}}
\Big(\sum_{i=1}^{n} \Big(y_{i} - \sum_{k=1}^{K} x_{i,k} \cdot q_{k} \Big)^{2}  
+ \lambda \cdot || \mathbf{Q} ||_{1} \Big),
\nonumber
\end{equation}

\vspace{-0.2cm}
\noindent
with data observations $y_{i}$ ($i=1,\ldots,n)$, deterministic explanatory variables $x_{i,k}$ and 
$\ell_{1}-$norm-regularization (penalization) parameter $\lambda \geq 0$.
This LASSO task is --- as e.g. discussed in Chapter 9 of Theodoridis 
\cite{BerBroStuGSI25:Theo20} ---
equivalent to finding the minimizer for the
\textit{basis pursuit denoising} problem
(cf. Donoho et al. \cite{BerBroStuGSI25:Donoho06a}, 
see also e.g. 
Cand{\`e}s et al. \cite{BerBroStuGSI25:Can06},
Lustig et al. \cite{BerBroStuGSI25:Lustig07},
Cand\`es \cite{BerBroStuGSI25:Can08a}, 
Cand\`es et al. \cite{BerBroStuGSI25:Can08b}, 
Goldstein \& Osher \cite{BerBroStuGSI25:Goldst09},
Zhang et al. \cite{BerBroStuGSI25:Zhang14},
Edgar et al. \cite{BerBroStuGSI25:Edgar19})
\vspace{-0.1cm}
\begin{eqnarray}
&& 
\min_{\mathbf{Q} \in \mathbf{\Omega}} \, 
||  \mathbf{Q} ||_{1}
\label{BerBroStuGSI25:BDP3 approx new}
\\[-0.2cm]
&&
\textrm{with} \quad
\mathbf{\Omega} := \Big\{ \mathbf{Q} \in \mathbb{R}^{K}: \, 
 \sum_{i=1}^{n} \Big(y_{i} - \sum_{k=1}^{K} x_{i,k} \cdot q_{k} \Big)^{2} \leq \varepsilon 
\Big\}
\nonumber
\end{eqnarray}

\enlargethispage{0.5cm}

\begin{figure}[h]  \vspace{-0.5cm}
\centering
\subfigure[ ]
{\includegraphics[width=0.78\textwidth]{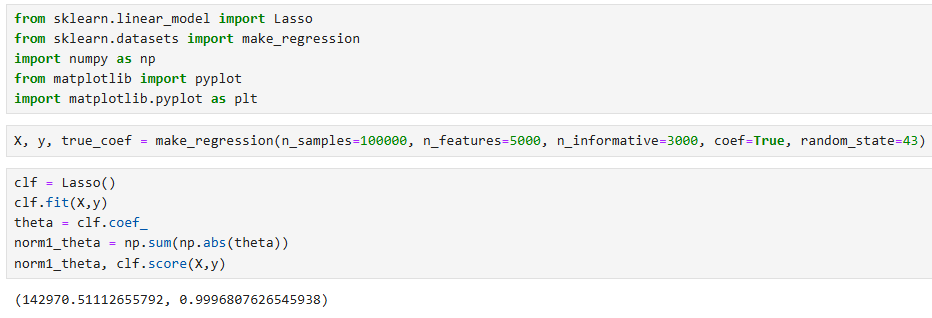}}
\hspace{0.05cm}
\subfigure[ ]
{\includegraphics[width=0.66\textwidth]{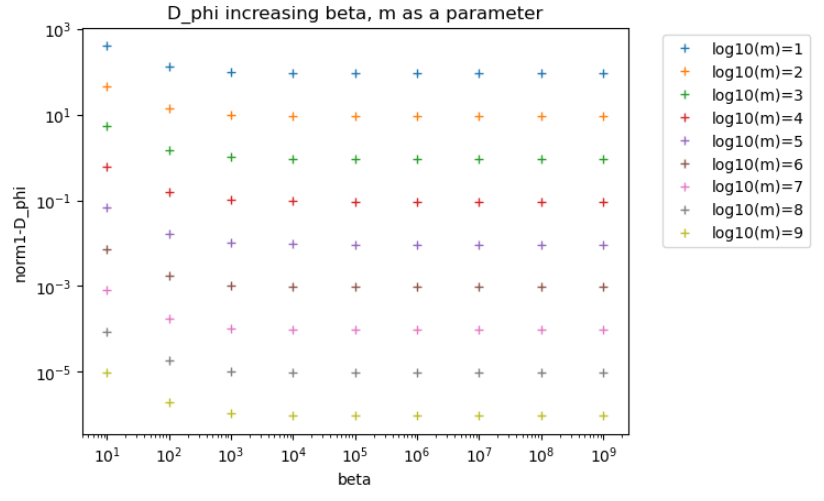}}
\hspace{0.05cm}
\subfigure[ ]
{\includegraphics[width=0.66\textwidth]{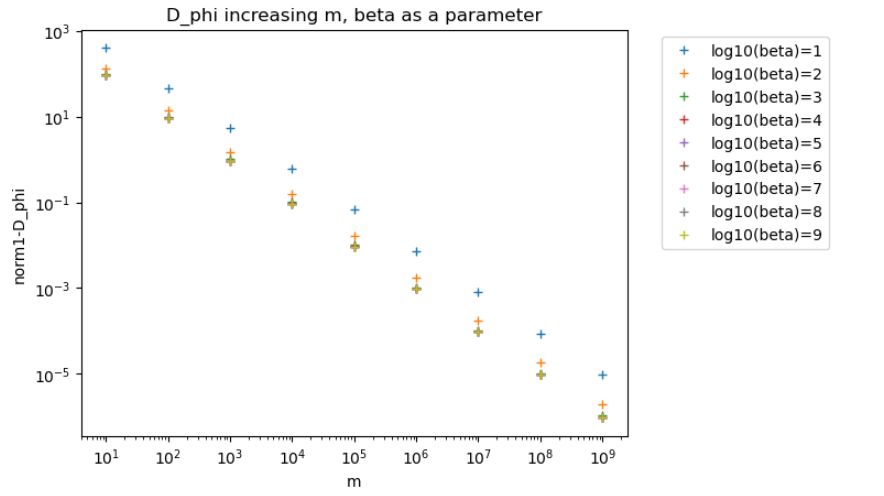}}
\hspace{0.05cm}
\subfigure[ ]
{\includegraphics[width=0.66\textwidth]{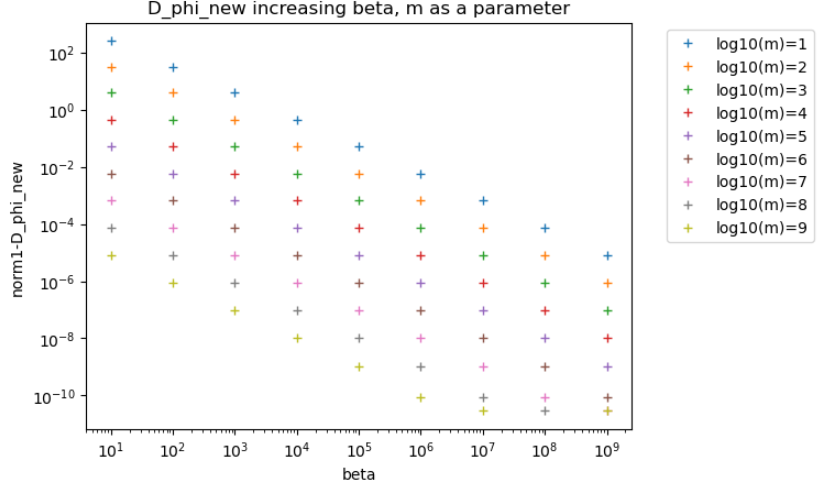}}
\vspace{-0.3cm}
\caption{\scriptsize \ }
\label{BerBroStuGSI25:fig1}
\end{figure}

\vspace{-0.1cm}
\noindent
for chosen fitting-quality parameter $\varepsilon>0$. In the light of this,
and in order to prepare for forthcoming studies dealing
with the (smooth bare-simulation-type, cf. 
\cite{BerBroStuGSI25:Bro23a},\cite{BerBroStuGSI25:Bro23b}) 
optimization of the corresponding smoothed version of
\eqref{BerBroStuGSI25:BDP3 approx new},
it is reasonable to compare $|| \widehat{\mathbf{Q}} ||_{1}$
with its smoother approximations 
$D_{\varphi_{\alpha,\beta,1/\beta}}\Big(\widehat{\mathbf{Q}},\frac{1}{m} \cdot \mathbf{1}\Big)$
and 
$D_{\varphi_{\alpha,\beta,1/\beta},\frac{1}{m} \cdot \mathbf{1},\mathbf{1}}^{new}(\widehat{\mathbf{Q}},\mathbf{0})$
for various different parameter constellations  $(\alpha,\beta, m)$
(cf. Proposition \ref{BerBroStuGSI25:prop.limit.ell1.BS3}(d)
and Proposition \ref{BerBroStuGSI25:prop.newdiv.ell1.BS3}(b)).
In the following, we analyse this for the LASSO-solution
generated by the Scikit-learn package by the code given in Figure~\ref{BerBroStuGSI25:fig1}(a),
where $K=5001$. 
Accordingly, we get $|| \widehat{\mathbf{Q}} ||_{1} = 142970.51$ 
(with a performance score of about $0.9997$).
In order to visually demonstrate the approximations
for some of the above-mentioned limits, let us always choose (say) $\alpha=1$.
Concerning \eqref{BerBroStuGSI25:fo.limit.ell1.2.BS3sec} 
with $\mathbf{P} = \frac{1}{m} \cdot \mathbf{1}$,
we plot --- on a logarithmic scale ---
in Figure~\ref{BerBroStuGSI25:fig1}(b) the function 
$\beta \ \mapsto \
|| \widehat{\mathbf{Q}} ||_{1} - 
D_{\varphi_{1,\beta,1/\beta}}(\widehat{\mathbf{Q}},\frac{1}{m} \cdot \mathbf{1})$
with
\vspace{-0.2cm}
\begin{equation}
\textstyle
D_{\varphi_{1,\beta,1/\beta}}(\widehat{\mathbf{Q}},\frac{1}{m} \cdot \mathbf{1})
=
\frac{1}{m \cdot \beta} \cdot  \sum\limits_{k=1}^{K} \negthinspace \negthinspace
\Big\{
\sqrt{1 + \beta^{2} \cdot (1 - m \cdot q_{k})^2} 
\, - \, 1 
+ \log\frac{2
}{
\sqrt{1 + \beta^{2} \cdot (1 - m \cdot q_{k})^2} \, + \, 1
}  \Big\}
\nonumber
\end{equation}
for increasingly large $\beta$, with several different values (in different colours) of large $m$
as a family parameter; in Figure~\ref{BerBroStuGSI25:fig1}(c), the roles of
$\beta$ and $m$ are switched. 
Concerning \eqref{BerBroStuGSI25:prop.newdiv.ell1.2.BS3sec} with $\mathbf{P} = \frac{1}{m} \cdot \mathbf{1}$,
$\mathbf{\sigma} = \mathbf{1}$ and $\mathbf{Q}^{\ast} = \mathbf{0}$,
we plot in Figure~\ref{BerBroStuGSI25:fig1}(d) the function 
$\beta \ \mapsto \ || \widehat{\mathbf{Q}} ||_{1} - 
D_{\varphi_{1,\beta,1/\beta},\frac{1}{m} \cdot \mathbf{1},\mathbf{1}}^{new}(\widehat{\mathbf{Q}},\mathbf{0})$
with
\begin{equation}
\textstyle
D_{\varphi_{1,\beta,1/\beta},\frac{1}{m} \cdot \mathbf{1},\mathbf{1}}^{new}(\widehat{\mathbf{Q}},\mathbf{0}) =
\frac{1}{m \cdot \beta} \cdot  \sum\limits_{k=1}^{K} \negthinspace 
\Big\{
\sqrt{1 + (m \cdot \beta \cdot q_{k})^2} 
\, - \, 1 
+ \log\frac{2
}{
\sqrt{1 + (m \cdot \beta \cdot q_{k})^2} \, + \, 1
}  \Big\}
\nonumber
\end{equation}
for increasingly large $\beta$, with several different values 
of large $m$ as a family parameter.

\vspace{-0.4cm}


\enlargethispage{0.4cm}

\section{Proofs}
\label{BerBroStuGSI25:sec:4}

\vspace{-0.3cm}
\noindent 
\textbf{Proof of Proposition \ref{BerBroStuGSI25:prop.limit.ell1.BS3}.} \ 
(a) 
For fixed $t \in \, ]-\infty,0[ \, \cup \, ]0,\infty[$, $\beta \in \, ]0,\infty[$ and $\widetilde{c} \in \, ]0,\infty[$,
one gets by \eqref{BerBroStuGSI25:fo.genLap3ba.BS3}
--- with the help of $y:=\frac{t}{\alpha}$ --- 
by De l'Hospital's rule
\vspace{-0.3cm}
\begin{eqnarray}
& & 
\textstyle
\hspace{-0.2cm} 
\lim_{\alpha \rightarrow 0_{+}} \varphi_{\alpha,\beta,\widetilde{c}}(t+1)  = 
\lim_{\alpha \rightarrow 0_{+}} 
\bigg( \widetilde{c} \cdot \alpha \cdot \Big\{
\sqrt{1 + \beta^{2} 
\cdot \frac{t^2}{\alpha^2}} \, - \, 1 
+ \log\frac{2 \cdot \Big(
\sqrt{1 + \beta^{2} 
\cdot \frac{t^2}{\alpha^2}} \, - \, 1
\Big)
}{
\beta^{2} \cdot \frac{t^2}{\alpha^2}
}  \Big\}
 \bigg)
\nonumber 
\\
& & 
\hspace{-0.2cm} 
= \ \widetilde{c} \cdot \Big\{
\beta \cdot |t| + 
\lim_{\alpha \rightarrow 0_{+}} 
\alpha \cdot \log\frac{2 \cdot \Big(
\sqrt{1 + \beta^{2} 
\cdot \frac{t^2}{\alpha^2}} \, - \, 1
\Big)
}{
\beta^{2} \cdot \frac{t^2}{\alpha^2}
}
\Big\}
\nonumber\\
& & 
\hspace{-0.2cm}
= \ \widetilde{c} \cdot \Big\{
\beta \cdot |t| + 
\lim_{|y| \rightarrow \infty} 
\frac{t}{y} \cdot 
\log\frac{2 \cdot \Big(
\sqrt{1 + \beta^{2} 
\cdot y^2} \, - \, 1
\Big)
}{
1 + \beta^{2} \cdot y^2 - 1
}
\Big\}
\label{BerBroStuGSI25:fo.proof.1.BS3}
\\
& & 
\hspace{-0.2cm}
=  \ \widetilde{c} \cdot \Big\{
\beta \cdot |t| + t \cdot
\lim_{|y| \rightarrow \infty} 
\frac{1}{y} \cdot 
\log\frac{2}{
\sqrt{1 + \beta^{2} \cdot y^2} + 1
}
\Big\}
\label{BerBroStuGSI25:fo.proof.2.BS3} \\
& &
\textstyle 
\hspace{-0.2cm}
=  \widetilde{c} \cdot \Big\{
\beta \cdot |t| + t \cdot \negthinspace
\lim_{|y| \rightarrow \infty} 
\frac{\sqrt{1 + \beta^{2} \cdot y^2}}{y}
\cdot 
\negthinspace \lim_{|y| \rightarrow \infty}
\frac{1}{\sqrt{1 + \beta^{2} \cdot y^2}} \cdot 
\log\frac{2}{
\sqrt{1 + \beta^{2} \cdot y^2} + 1
}
\Big\}
 = 
\widetilde{c} \cdot \beta \cdot |t| . \hspace{0.6cm} {\ }
\label{BerBroStuGSI25:fo.proof.3.BS3}
\end{eqnarray}

\vspace{-0.2cm}
\noindent
Moreover, for $t=0$ one has $\varphi_{\alpha,\beta,\widetilde{c}}(t+1) = 
\varphi_{\alpha,\beta,\widetilde{c}}(1) = 0$ even for all $\alpha \in \, ]0,\infty[$.\\
(b) This works analogously to the proof of (a).\\
(c) By means of (a) we get for all $\beta \in \, ]0,\infty[$, $\widetilde{c} \in \, ]0,\infty[$,
$\mathbf{Q} \in \mathbb{R}^{K}$ and $\mathbf{P} \in \mathbb{R}_{>0}^{K}$
\vspace{-0.2cm}
\begin{eqnarray}
& & \lim_{\alpha \rightarrow 0_{+}}
D_{\varphi_{\alpha,\beta,\widetilde{c}}}(\mathbf{Q},\mathbf{P}) \ = \ 
\sum\limits_{k=1}^{K} 
p_{k} \cdot 
\lim_{\alpha \rightarrow 0_{+}}
\varphi_{\alpha,\beta,\widetilde{c}}\Big(\frac{q_{k}}{p_{k}}\Big)
\ = \ 
\widetilde{c} \cdot \beta \cdot \sum_{k=1}^{K} 
p_{k} \cdot \Big|\frac{q_{k}}{p_{k}} - 1\Big| . 
\nonumber
\end{eqnarray}

\vspace{-0.2cm}
\noindent
(d) This works analogously to the proof of (c).\\
(e) Let us arbitrarily fix $\alpha \in \, ]0,\infty[$, $\beta \in \, ]0,\infty[$
and $\widetilde{c} \in\, ]0,\infty[$.
For all
$p \in \, ]0,\infty[$
and $q \in \, ]0,\infty[$
($k=1,\ldots,K$)
one can derive --- by setting $x: = \frac{q}{p} >0$ 
and $y:= \frac{x-1}{\alpha}$ ---
by De l'Hospital's rule
\vspace{-0.2cm}
\begin{eqnarray}
& & 
\hspace{-0.2cm}
\lim_{p \rightarrow 0_{+}} 
p \cdot \varphi_{\alpha,\beta,\widetilde{c}}\Big(\frac{q}{p}\Big)  =  
\lim_{p \rightarrow 0_{+}} 
q \cdot \frac{1}{\frac{q}{p}} \cdot \varphi_{\alpha,\beta,\widetilde{c}}\Big(\frac{q}{p}\Big) 
= q \cdot \lim_{x \rightarrow\infty} 
\frac{1}{x} \cdot \varphi_{\alpha,\beta,\widetilde{c}}(x)
\nonumber \\
& &
\hspace{-0.2cm}
\textstyle 
= 
q \cdot \lim_{x \rightarrow \infty} 
\bigg( \widetilde{c} \cdot \frac{x-1}{x} \cdot \frac{\alpha}{x-1} \cdot \Big\{
\sqrt{1 + \beta^{2} 
\cdot \frac{(x-1)^2}{\alpha^2}} \, - \, 1 
+ \log\frac{2 \cdot \Big(
\sqrt{1 + \beta^{2} 
\cdot \frac{(x-1)^2}{\alpha^2}} \, - \, 1
\Big)
}{
\beta^{2} \cdot \frac{(x-1)^2}{\alpha^2}
}  \Big\}
 \bigg)
\nonumber \\
& & 
\hspace{-0.2cm}
= \ q \cdot \widetilde{c} \cdot \Big\{
\beta + 
\lim_{y \rightarrow \infty} 
\frac{1}{y} \cdot 
\log\frac{2 \cdot \Big(
\sqrt{1 + \beta^{2} 
\cdot y^2} \, - \, 1
\Big)
}{
1 + \beta^{2} \cdot y^2 - 1
}
\Big\}
\ = \ \ q \cdot \widetilde{c} \cdot \beta,
\nonumber
\end{eqnarray}
where the last equality follows as in \eqref{BerBroStuGSI25:fo.proof.1.BS3},
\eqref{BerBroStuGSI25:fo.proof.2.BS3}
and \eqref{BerBroStuGSI25:fo.proof.3.BS3} above.
Analogously, for all
$p \in \, ]0,\infty[$
and $q \in \, ]-\infty,0[$
one can derive --- by setting $x: = - \frac{q}{p} >0$ 
and $y:= \frac{x+1}{\alpha}$ ---
by De l'Hospital's rule
\vspace{-0.2cm}
\begin{eqnarray}
& & \hspace{-0.2cm}
\lim_{p \rightarrow 0_{+}} 
p \cdot \varphi_{\alpha,\beta,\widetilde{c}}\Big(\frac{q}{p}\Big) \ = \ 
\lim_{p \rightarrow 0_{+}} 
q \cdot \frac{1}{\frac{q}{p}} \cdot \varphi_{\alpha,\beta,\widetilde{c}}\Big(\frac{q}{p}\Big) 
\ = \ - \, q \cdot \lim_{x \rightarrow \infty} 
\frac{1}{x} \cdot \varphi_{\alpha,\beta,\widetilde{c}}(-x)
\nonumber \\
& & \hspace{-0.2cm}
\textstyle
\ = \ 
- \, q \cdot \lim_{x \rightarrow \infty} 
\bigg( \widetilde{c} \cdot \frac{x+1}{x} \cdot \frac{\alpha}{x+1} \cdot \Big\{
\sqrt{1 + \beta^{2} 
\cdot \frac{(x+1)^2}{\alpha^2}} \, - \, 1 
+ \log\frac{2 \cdot \Big(
\sqrt{1 + \beta^{2} 
\cdot \frac{(x+1)^2}{\alpha^2}} \, - \, 1
\Big)
}{
\beta^{2} \cdot \frac{(x+1)^2}{\alpha^2}
}  \Big\}
 \bigg)
\nonumber \\
& & \hspace{-0.2cm}
= \ - \, q \cdot \widetilde{c} \cdot \Big\{
\beta + 
\lim_{y \rightarrow \infty} 
\frac{1}{y} \cdot 
\log\frac{2 \cdot \Big(
\sqrt{1 + \beta^{2} 
\cdot y^2} \, - \, 1
\Big)
}{
1 + \beta^{2} \cdot y^2 - 1
}
\Big\}
\ = \ - \, q \cdot \widetilde{c} \cdot \beta.
\nonumber
\end{eqnarray}

\vspace{-0.2cm}
\noindent
Moreover, for all $p \in \, ]0,\infty[$
and $q = 0$ one has $\lim_{p \rightarrow 0_{+}} 
p \cdot \varphi_{\alpha,\beta,\widetilde{c}}\Big(\frac{q}{p}\Big)
= \lim_{p \rightarrow 0_{+}} 
p \cdot \varphi_{\alpha,\beta,\widetilde{c}}(0) = 0 = q \cdot \widetilde{c} \cdot \beta$
since $\varphi_{\alpha,\beta,\widetilde{c}}(0) \in \, ]0,\infty[$.
Summing up, we have shown for all
$p \in \, ]0,\infty[$
and $q \in \, ]-\infty,\infty[$ that 
$\lim_{p \rightarrow 0_{+}} 
p \cdot \varphi_{\alpha,\beta,\widetilde{c}}\Big(\frac{q}{p}\Big) =  
\widetilde{c} \cdot \beta \cdot |q|$.
From this and the notation $\mathbf{P}_{m} := (p_{m,1}, \ldots, p_{m,K})$, 
\eqref{BerBroStuGSI25:fo.limit.ell1.2.BS3a} follows immediately from
\vspace{-0.3cm}
\begin{eqnarray}
& & \lim_{m \rightarrow \infty}
D_{\varphi_{\alpha,\beta,\widetilde{c}}}(\mathbf{Q},\mathbf{P}_{m}) \ = \ 
\lim_{m \rightarrow \infty}
\sum\limits_{k=1}^{K} 
p_{m,k} \cdot 
\varphi_{\alpha,\beta,\widetilde{c}}\Big(\frac{q_{k}}{p_{m,k}}\Big)
\ = \ 
\widetilde{c} \cdot \beta \cdot \sum_{k=1}^{K} |q_{k}|, 
\nonumber
\end{eqnarray}

\vspace{-0.2cm}
\noindent
where we have employed that $p_{m,k} \stackrel{m\rightarrow \infty}{\longrightarrow} 0$
for all $k=1,\ldots,K$. \hspace{0.5cm}  $\blacksquare$ 

\vspace{0.2cm}
\noindent
\textbf{Proof of Proposition \ref{BerBroStuGSI25:prop.newdiv.ell1.BS3}.} \ 
(a) By means of Proposition \ref{BerBroStuGSI25:prop.limit.ell1.BS3}(a) we get for all $\beta \in \, ]0,\infty[$, $\widetilde{c} \in \, ]0,\infty[$,
$\mathbf{Q} \in \mathbb{R}^{K}$, $\mathbf{Q}^{\ast} \in \mathbb{R}^{K}$, $\mathbf{P} \in \mathbb{R}_{>0}^{K}$
and $\mathbf{\sigma} \in \mathbb{R}_{>0}^{K}$
\vspace{-0.2cm}
\begin{equation}
\lim_{\alpha \rightarrow 0_{+}}
D_{\varphi_{\alpha,\beta,\widetilde{c}},\mathbf{P},\mathbf{\sigma}}^{new}(\mathbf{Q},\mathbf{Q}^{\ast})
= \sum\limits_{k=1}^{K} p_{k} \cdot
\lim_{\alpha \rightarrow 0_{+}}
\varphi_{\alpha,\beta,\widetilde{c}} \negthinspace 
\left( \frac{q_{k}-q_{k}^{\ast}}{p_{k} \cdot \sigma_{k}} + 1 \right)
 = \widetilde{c} \cdot \beta \cdot \sum_{k=1}^{K} 
p_{k} \cdot \Big|\frac{q_{k}-q_{k}^{\ast}}{p_{k} \cdot \sigma_{k}}\Big| .
\nonumber
\end{equation}

\vspace{-0.2cm}
\noindent
(b) This works analogously to the proof of (a), by employing 
Proposition \ref{BerBroStuGSI25:prop.limit.ell1.BS3}(b).\\
(c) Let us arbitrarily fix $\alpha \in \, ]0,\infty[$, $\beta \in \, ]0,\infty[$
and $\widetilde{c} \in\, ]0,\infty[$.
For all
$p \in \, ]0,\infty[$
and $\breve{q} \in \, ]0,\infty[$
one can derive --- by setting $x: = \frac{\breve{q}}{p} >0$ 
and $y:= \frac{x}{\alpha} >0$ ---
by De l'Hospital's rule
\vspace{-0.2cm}
\begin{eqnarray}
& & \hspace{-0.2cm}
\lim_{p \rightarrow 0_{+}} 
p \cdot \varphi_{\alpha,\beta,\widetilde{c}}\Big(\frac{\breve{q}}{p}+1\Big)  = 
\lim_{p \rightarrow 0_{+}} 
\breve{q} \cdot \frac{1}{\frac{\breve{q}}{p}} \cdot \varphi_{\alpha,\beta,\widetilde{c}}\Big(\frac{\breve{q}}{p}+1\Big) 
= \breve{q} \cdot \lim_{x \rightarrow\infty} 
\frac{1}{x} \cdot \varphi_{\alpha,\beta,\widetilde{c}}(x+1)
\nonumber \\
& & \hspace{-0.2cm}
= \ 
\breve{q} \cdot \lim_{x \rightarrow \infty} 
\bigg( \widetilde{c} \cdot \frac{\alpha}{x} \cdot \Big\{
\sqrt{1 + \beta^{2} 
\cdot \frac{x^2}{\alpha^2}} \, - \, 1 
+ \log\frac{2 \cdot \Big(
\sqrt{1 + \beta^{2} 
\cdot \frac{x^2}{\alpha^2}} \, - \, 1
\Big)
}{
\beta^{2} \cdot \frac{x^2}{\alpha^2}
}  \Big\}
 \bigg)
\nonumber \\
& & \hspace{-0.2cm}
= \ \breve{q} \cdot \widetilde{c} \cdot \Big\{
\beta + 
\lim_{y \rightarrow \infty} 
\frac{1}{y} \cdot 
\log\frac{2 \cdot \Big(
\sqrt{1 + \beta^{2} 
\cdot y^2} \, - \, 1
\Big)
}{
1 + \beta^{2} \cdot y^2 - 1
}
\Big\}
\ = \ \breve{q} \cdot \widetilde{c} \cdot \beta
\ = \ |\breve{q}| \cdot \widetilde{c} \cdot \beta,
\nonumber
\end{eqnarray}

\vspace{-0.2cm}
\noindent
where the last equality follows as in
\eqref{BerBroStuGSI25:fo.proof.1.BS3},
\eqref{BerBroStuGSI25:fo.proof.2.BS3}
and \eqref{BerBroStuGSI25:fo.proof.3.BS3} 
above. Analogously, for all
$p \in \, ]0,\infty[$
and $\breve{q}\in \, ]-\infty,0[$
one can derive --- by setting $x: = - \frac{\breve{q}}{p} >0$ 
and $y:= \frac{x}{\alpha}>0$ ---
by De l'Hospital's rule
\vspace{-0.2cm}
\begin{eqnarray}
& & 
\lim_{p \rightarrow 0_{+}} 
p \cdot \varphi_{\alpha,\beta,\widetilde{c}}\Big(\frac{\breve{q}}{p}+1\Big) \ = \ 
\lim_{p \rightarrow 0_{+}} 
\breve{q} \cdot \frac{1}{\frac{\breve{q}}{p}} \cdot \varphi_{\alpha,\beta,\widetilde{c}}\Big(\frac{\breve{q}}{p}+1\Big) 
\ = \ - \, \breve{q} \cdot \lim_{x \rightarrow \infty} 
\frac{1}{x} \cdot \varphi_{\alpha,\beta,\widetilde{c}}(-x+1)
\nonumber \\
& &
\ = \ 
- \, \breve{q} \cdot \lim_{x \rightarrow \infty} 
\bigg( \widetilde{c} \cdot \frac{\alpha}{x} \cdot \Big\{
\sqrt{1 + \beta^{2} 
\cdot \frac{x^2}{\alpha^2}} \, - \, 1 
+ \log\frac{2 \cdot \Big(
\sqrt{1 + \beta^{2} 
\cdot \frac{x^2}{\alpha^2}} \, - \, 1
\Big)
}{
\beta^{2} \cdot \frac{x^2}{\alpha^2}
}  \Big\}
 \bigg)
\nonumber \\
& & 
= \ - \, \breve{q} \cdot \widetilde{c} \cdot \Big\{
\beta + 
\lim_{y \rightarrow \infty} 
\frac{1}{y} \cdot 
\log\frac{2 \cdot \Big(
\sqrt{1 + \beta^{2} 
\cdot y^2} \, - \, 1
\Big)
}{
1 + \beta^{2} \cdot y^2 - 1
}
\Big\}
\ = \ - \, \breve{q} \cdot \widetilde{c} \cdot \beta
\ = \ |\breve{q}| \cdot \widetilde{c} \cdot \beta.
\nonumber
\end{eqnarray}

\vspace{-0.2cm}
\noindent
Moreover, for all $p \in \, ]0,\infty[$
and $\breve{q} = 0$ one has 
$p \cdot \varphi_{\alpha,\beta,\widetilde{c}}\Big(\frac{\breve{q}}{p}+1\Big)
= p \cdot \varphi_{\alpha,\beta,\widetilde{c}}(1) = 0 = |\breve{q}| \cdot \widetilde{c} \cdot \beta$.
Summing up, we have shown for all
$p \in \, ]0,\infty[$
and $\breve{q} \in \, ]-\infty,\infty[$ that 
$\lim_{p \rightarrow 0_{+}} 
p \cdot \varphi_{\alpha,\beta,\widetilde{c}}\Big(\frac{\breve{q}}{p}+1\Big) =  
\widetilde{c} \cdot \beta \cdot |\breve{q}|$.
From this and the notation $\mathbf{P}_{m} := (p_{m,1}, \ldots, p_{m,K})$, 
the desired limit relation \eqref{BerBroStuGSI25:prop.newdiv.ell1.2.BS3a} follows immediately from
\vspace{-0.2cm}
\begin{equation}
\lim_{m \rightarrow \infty}
D_{\varphi_{\alpha,\beta,\widetilde{c}},\mathbf{P}_{m},\mathbf{\sigma}}^{new}(\mathbf{Q},\mathbf{Q}^{\ast})
= \lim_{m \rightarrow \infty}
\sum\limits_{k=1}^{K} 
p_{m,k} \cdot 
\varphi_{\alpha,\beta,\widetilde{c}}\Big(\frac{q_{k} - q_{k}^{\ast}}{\sigma_{k} \cdot p_{m,k}} +1\Big)
= \widetilde{c} \cdot \beta \cdot \sum_{k=1}^{K} \frac{|q_{k} - q_{k}^{\ast}|}{\sigma_{k}},
\nonumber
\end{equation}

\vspace{-0.2cm}
\noindent
where we have employed that $p_{m,k} \stackrel{m\rightarrow \infty}{\longrightarrow} 0$
for all $k=1,\ldots,K$. \hspace{0.5cm}  $\blacksquare$

\vspace{-0.2cm}


\subsubsection{Acknowledgements}  
W. Stummer is grateful to the Sorbonne Universit\'{e} 
Paris for its multiple partial financial support and especially to the LPSM 
for its multiple great hospitality.

\enlargethispage{0.5cm}

\vspace{-0.2cm}


\newpage

\end{document}